\documentclass[12pt]{amsart}
\usepackage{amscd,amssymb}
\usepackage[graph,frame,poly,arc]{xy}  
\usepackage[plainpages,backref,urlcolor=blue]{hyperref}

\theoremstyle{plain}
\newtheorem{thm}[subsection]{Theorem}
\newtheorem{lem}[subsection]{Lemma}
\newtheorem{prop}[subsection]{Proposition}
\newtheorem{cor}[subsection]{Corollary}

\theoremstyle{definition}
\newtheorem{rk}[subsection]{Remark}
\newtheorem{definition}[subsection]{Definition}
\newtheorem{ex}[subsection]{Example}

\numberwithin{equation}{section}
\newcommand{\OO}{{\mathcal O}}

\newcommand{\M}{{\mathcal M}}

\newcommand{\A}{{\mathcal A}}
\newcommand{\wJ}{\widehat{{J}}}

\newcommand{\owari}{\hfill$\square$}

\newcommand{\B}{{\mathcal B}}

\newcommand{\al}{{\alpha}}
\newcommand{\be}{{\beta}}

\newcommand{\Z}{\mathbb{Z}}

\newcommand{\C}{\mathbb{C}}
\newcommand{\K}{\mathbb{K}}

\newcommand{\PP}{\mathbb{P}}

\DeclareMathOperator{\Der}{Der}

\DeclareMathOperator{\nexp}{nexp}

\DeclareMathOperator{\im}{im}
\DeclareMathOperator{\coker}{coker}

\begin{document}

\title [Splitting types of bundles of logarithmic vector fields]
{On the splitting types of  bundles of logarithmic vector fields along plane curves}

\author[Takuro Abe and Alexandru Dimca]
{Takuro Abe\address{
Institute of Mathematics for Industry, 
Kyushu University,
Fukuoka 819-0395, Japan}
\email{abe@imi.kyushu-u.ac.jp} 
and Alexandru Dimca 
\address{Universit\'e C\^ ote d'Azur, CNRS, LJAD, France}
\email {dimca@unice.fr}
}

\maketitle

\begin{abstract}
We give a formula relating the total Tjurina number and 
the generic splitting type of the bundle of logarithmic vector fields associated to a reduced plane 
curve. By using it, we give a characterization of nearly free curves in terms of 
splitting types.  Several applications to free and nearly free arrangements of lines are also given,
in particular a proof of a form of Terao's Conjecture for arrangements having a line with at most 4 intersection points.
\end{abstract}

\section{Introduction} \label{s0}

Let $C:f=0$ be a reduced curve of degree $d$ in $X=\PP^2$, $S=\C[x,y,z]$, and $AR(f)$ the graded $S$-module of Jacobian syzygies of $f$ as in \cite{DS14}, see equation \eqref{eqAR} below. Note that 
$AR(f)$ is isomorphic to the logarithmic derivation module $D_0(C)$ of the 
curve $C$ defined by $D_0(C):=\{\theta \in \mbox{Der} S \mid 
\theta(f)=0\}$.  Let  $mdr(f):=\min \{k \mid 
AR(f)_k \neq (0)\}$. In this paper we assume that 
$$
mdr(f) \ge 1,
$$
unless otherwise specified. The case $mdr(f)=0$ corresponds to the rather trivial case when $C$ is a collection of $d$ lines passing through one common point.
Let $E_C$ be the locally free sheaf on $X$ corresponding to the graded module $AR(f)$, and recall that
\begin{equation} \label{equa1} 
E_C=T\langle C \rangle (-1),
\end{equation}
where $T\langle C \rangle $ is the sheaf of logarithmic vector fields along $C$ as considered for instance in \cite{DS14}.
For a line $L$ in $X$, the pair of integers $( d_1^L, d_2^L)$, with $ d_1^L \leq d_2^L$, such that $ E_C|_L \simeq \OO_L(-d_1^L) \oplus \OO_L(-d_2^L)$ is called the splitting type of $E_C$ along $L$, see for instance \cite{FV1,OSS}. For a generic line $L_0$, the corresponding splitting type $( d_1^{L_0}, d_2^{L_0})$ is constant. For any line $L$ in $X$, we set
$$I(C,L)=(d-1)^2-d_1^{L} d_2^{L}.$$ 

The algebraic structure of the graded $S$-module $AR(f)$ is related to the 
singularities of $C$, e.g. the invariants like Milnor numbers and Tjurina numbers. From this viewpoint, when the $S$-module $AR(f)$ is free, 
which can be considered as the simplest case, then the corresponding curve is called free, a notion going back to K. Saito \cite{S}.
When the minimal resolution of the graded $S$-module $AR(f)$ is slightly more complicated, we get the nearly free curves considered in \cite{DSt15}. 
See \S 2 for details. 
Recall the definition of the global Tjurina number 
$$\tau(C)= \sum _{p \in C} \tau(C,p)$$ 
of the curve $C$. 
Also, let $N(f)= \wJ_f/J_f$, with $J_f$ the Jacobian ideal of $f$ in $S$, spanned by the partial derivatives $f_x,f_y$, $f_z$ of $f$, and $\wJ_f$ the saturation of the ideal $J_f$ with respect to the maximal ideal 
${\bf m}=(x,y,z)$ in $S$. The quotient module $N(f)=H^0_{\bf m}(S/J_f)$ plays a key role in this theory.
Indeed, let $\nu(C)= \dim N(f)_{[T/2]}$, where $[ \ \ ]$ denotes integral part and $T=3(d-2)$.
It is known that the curve $C:f=0$ is free (resp. nearly free) if and only if $\nu(C)=0$ (resp. $\nu(C)=1$), see  \cite{D2, DPop, DSt15}.
The above key invariants of distinct origins are related in the first main result of this paper.
\begin{thm} \label{thm1} 
 With the above notation,  for any line $L$ in $\PP^2$, and any generic line $L_0$ in $\PP^2$,  the following hold.
 \begin{enumerate}
 
 \item $\max(mdr(f)-\nu(C),0) \leq d_1^{L} \leq  d_1^{L_0} \leq  \min(mdr(f),[(d-1)/2]).$
 
\item  $I(C,L) \geq I(C,L_0)=\tau(C) + \nu(C).$

 \end{enumerate}
In particular, the reduced curve $C:f=0$ in $\PP^2$ is free (resp. nearly free) if and only if $I(C,L_0)=\tau(C)$ (resp. $I(C,L_0)=\tau(C)+1$).

\end{thm}

\begin{cor}\label{mainAG}
Let 
$c_{E_C}(t)=1+c_1(E_C)t+c_2(E_C) t^2\in \Z[t]$ be the Chern polynomial of the vector bundle $E_C$. Then the curve $C$ is 
free (resp. nearly free) if and only if there is a line 
$L \subset \PP^2$ such that 
$c_2(E_C)-d_1^Ld_2^L=0$
(resp. $c_2(E_C)-d_1^Ld_2^L=1$).

\end{cor}

The free case of Corollary \ref{mainAG} 
is due to Yoshinaga in \cite{Y3}, see Theorem \ref{Y} below. 
We give another proof for this case in terms of Tjurina numbers. 
Indeed, the proof of Theorem \ref{thm1}  implies
$$c_2(E_C)-d_1^{L_0} d_2^{L_0}=(d-1)^2-\tau(C)-d_1^{L_0} d_2^{L_0}=\nu(C).$$
This relation also yields the following analog of a result in \cite{duPCTC}.
\begin{cor}\label{corCTC}
For any reduced curve $C:f=0$ in $\PP^2$ and any line $L$ in $\PP^2$ we have
$$\tau(C) \leq (d-1)^2-d_1^{L} d_2^{L}.$$
Moreover, equality holds for a line $L$ if and only if the curve $C$ is free, and then it holds for any line $L$.
\end{cor}

The second main result of our paper is the following.

\begin{thm}\label{mdr}
Let $C:f=0$ be a reduced curve of degree $d$ in projective plane $\PP^2$.
Then the following properties are equivalent.
\begin{enumerate}
\item  The Chern polynomial $c_{E_C}(t)=1+c_1(E_C)t+c_2(E_C) t^2\in \Z[t]$ of the vector bundle $E_C$ has real roots.

\item $\tau(C)\geq \frac{3}{4}(d-1)^2.$

\item $d_1^{L_0} \leq \frac{d-1}{2}- \sqrt {\nu(C)}$ for a generic line $L_0$.

\end{enumerate}
Moreover, these properties  imply that
$mdr(f)$ coincides with $d_1^{L_0}$ for a generic line $L_0$, and they are implied by either
$mdr(f) <\frac{d}{4}$     or    
$$
 mdr(f) \leq \frac{d-1}{2}- \sqrt {\nu(C)}.$$
\end{thm}

The organization of this paper is as follows. In \S 2 we recall several 
definitions and results necessary for the proof of the main results. In \S 3 we prove Theorem \ref{thm1} 
and  Corollary \ref{mainAG}.  In \S 4, we  prove Theorem \ref{mdr}. 
In \S 5 and \S 6 we apply these results to the case of a line arrangement $\A:f=0$ in 
$\PP^2$. A sample of the results we get in this case is the following special case of Terao's Conjecture.

\begin{thm}
Let $\A$ be a line arrangement in $\PP^2$, such that
\begin{enumerate}
\item  there is a line $H \in \A$ containing at most 4 intersection points of $\A$; 

\item this line  $H$ does not contain an intersection point of multiplicity  $ \geq |\A|/2$. 
\end{enumerate}
Then the fact that $\nu(\A) \leq 1$, i.e. the fact that $\A$ is either free or nearly free, 
depends only on the characteristic polynomial
$\chi(\A;t)$. 
In addition, when the characteristic polynomial
$\chi(\A;t)$ is not a perfect square, it determines precisely whether $\A$ is free or nearly free.
\label{456}
\end{thm}


\medskip

\noindent
\textit{Acknowledgements}. 
A part of the argument in the proof of Theorems \ref{addition3} and 
\ref{addition4} is due to 
the first author's joint work with Max Wakefield. The authors 
are really grateful to him for his letting 
the argument to be used in this paper. 
The first author is 
partially supported by JSPS Grant-in-Aid for Scientific Research (B) 16H03924, and 
Grant-in-Aid for Exploratory Research 16K13744.

\section{Preliminaries}

First let us recall the definition of free and nearly free curves. In this 
paper $C$ is a reduced plane curve in $\PP^2$, defined by $f=0$, 
where $f$ is a degree $d$ homogeneous polynomial. For the coordinate ring 
$S=\C[x,y,z]$ and a graded $S$-module $M$, let $M_k$ denote the homogeneous 
degree $k$-part of $M$. For $g \in S$, let $g_x,g_y,g_z$ denote the partial 
derivative of $g$ by $x,y,z$. Then the graded $S$-module $AR(f)=AR(C) 
\subset S^{\oplus 3}$ of 
\textit{all relations} is defined by 
\begin{equation} \label{eqAR} 
AR(f)_k:=\{(a,b,c) \in S^{\oplus 3}_k \mid a f_x+b f_y+c f_z=0\}.
\end{equation}
The module $AR(f)$ is  isomorphic to the logarithmic derivations killing $f$, hence we sometimes identify the two types of objects as follows:
$$
AR(f) \ni (a,b,c) \mapsto a \partial_x + b \partial_y+ c \partial_z \in D_0(C).
$$
Its sheafification $E_C:=
\widetilde{AR(f)}$ is a rank two vector bundle on $\PP^2$, see 
\cite{S, Se} for details. In particular we have the following.

\begin{prop}[Equation (3.1), \cite{DS14}]
For a coherent sheaf $F$ on $\PP^2$, consider the graded $S$-module
$\Gamma_*(F):=
\oplus_{k \in \Z} H^0(\PP^2,F(k)).$ 
Then 
$
\Gamma_*(E_C)=AR(f).
$
\label{ID}
\end{prop}

\begin{definition}[\cite{DSt15}]

(1)\,\,
A curve $C$ is \textit{free} if the graded $S$-module
$AR(f)$ is free, say with  a  basis $\varphi_1,\varphi_2$.  If
$\deg \varphi_i=d_i\ (i=1,2)$,
the multiset of integers $(d_1,d_2)$ is called the 
\textit{exponents} of a free curve $C$, and denoted by $\exp(C)=(d_1,d_2)$.

(2)\,\,
A curve $C$ is \textit{nearly free} if 
$N(f) \neq 0$ and $\dim N(f)_k \le1$ for any $k$.
\label{nearlyfree}
\end{definition}

By \cite{DSt15}, the near freeness coincides with the following.

\begin{prop}[\cite{DSt15}]
$C$ is nearly free if and only if the graded $S$-module
$AR(f)$ has a minimal generator system of syzygies $\theta, \varphi_1,\varphi_2$, such that 
$$\deg \theta=
d_1 \le d_2:=\deg \varphi_1=\deg \varphi_2$$ 
with a relation 
$$
h \theta +\beta_1 \varphi_1+\beta_2 \varphi_2=0,
$$
for $h \in S$ and linear forms $\beta_1,\beta_2$. The multiset $(d_1,d_2)$ is called the
\textit{exponents} of a nearly free curve $C$, 
and denoted by $\nexp(C)=(d_1,d_2)$.
\label{nearlyfree2}
\end{prop}

Hence in terms of sheaves, for a nearly free curve $C$, 
the bundle $E_C$ has a minimal resolution of the form
$$
0 \rightarrow 
\mathcal{O}(-d_2-1)
\rightarrow 
\mathcal{O}(-d_1)\oplus \mathcal{O}(-d_2)^{\oplus 2} 
\rightarrow E_C \rightarrow 0.
$$
The following statement is immediate.

\begin{prop}
For a nearly free curve with $\nexp(C)=(d_1,d_2)$, it holds that 
$$
c_t(E_C)=1-(d_1+d_2-1)t+(d_1(d_2-1)+1)t^2.
$$
\label{Chern}
\end{prop}

Recall also the following characterization of nearly free curves.

\begin{prop}[Proposition 3.8, \cite{DSt15}]
$C$ is nearly free if and only if 
$$\nu(C)=1.$$
\label{dimNf}
\end{prop}

For the proof of the main results, we need the following. 
Let $\alpha_L$ be the defining equation of the line $L$. Then one has an exact sequence
$$0 \to \OO_X(-1) \stackrel{\cdot \alpha_L}{\to} \OO_X \to \OO_L \to 0,$$
where the first non-trivial morphism is induced by multiplication by the linear form $\alpha_L$. Let $k$ be an integer and tensor the above exact sequence by the vector bundle $E_C(k)$. We get
$$0 \to E_C(k-1) \stackrel{\cdot \alpha_L}{\to} E_C(k) \to E_C(k)|_L \to 0,$$
with $E_C(k)|_L \simeq \OO_L(k-d_1^L) \oplus \OO_L(k-d_2^L)$, since we assume as in the Introduction that $E_C|_L \simeq \OO_L(-d_1^L) \oplus \OO_L(-d_2^L)$. 
Then we have the following.

\begin{prop}
The long exact sequence of cohomology groups of the short exact sequence above starts as follows:
\begin{eqnarray} \label{equa2} 
0 &\to& AR(f)_{k-1} \stackrel{\cdot \alpha_L}{\to} AR(f)_{k} \stackrel{\pi_L}{\to} H^0(L,\OO_L(k-d_1^L) \oplus \OO_L(k-d_2^L)) \\
&\to& N(f)_{k+d-2} \stackrel{\cdot \alpha_L}{\to}
N(f)_{k+d-1} \to \cdots. \nonumber
\end{eqnarray}
Moreover, for $k=-1$, the corresponding morphism $N(f)_{d-3} \stackrel{\cdot \alpha_L}{\to}
N(f)_{d-2}$ is injective and $d_1^L\geq 0$  for any line $L$. 
\label{key1}
\end{prop}

\proof
This is exactly as in the proof of \cite[Theorem 5.7]{DS14}.  The key point is the identification $H^1(\PP^2, E_C(k))=N(f)_{k+d-1}$, valid for any integer $k$,  for which we refer 
to \cite[Proposition 2.1]{Se}. For the last claim, note that $N(f)_{d-3} \subset S_{d-3}$ and
$N(f)_{d-2} \subset S_{d-2}$, as the Jacobian ideal is generated in degree $d-1$.
\endproof

The following result is often used to investigate the structure of $AR(f)$. 

\begin{thm}[\cite{Y3}, Theorem 1.45]
Let $F$ be a rank two vector bundle on $\PP^2$ and $L \subset \PP^2$ 
a line. Then 
$$
\dim \coker (\pi_L: \Gamma_*(F) \rightarrow \Gamma_*(F|_L))=c_2(F) - d_1^Ld_2^L\ge 0,
$$ and 
the equality holds if and only if $F \simeq 
\mathcal{O}(-d_1^L) \oplus \mathcal{O}(-d_2^L)$.
\label{Y2}
\end{thm}

For a rank two vector bundle $E$ on $\PP^2$, consider the function
\begin{equation} \label{aE}
a_E:(\PP^2)^* \rightarrow \Z^2,  \text{ defined by }a_E(L):=
(d_1^L,d_2^L),
\end{equation}
where we take $d_1^L \leq d_2^L$.
 We say that $E$ is \textit{uniform} if the function $a_E$ is constant. 
The following classification of the uniform 2-bundles on $\PP^2$ is often used.

\begin{thm}[e.g., \S 2.2, Theorem 2.2.2, \cite{OSS}]
A rank two uniform vector bundle on $\PP^2$ is either 
(a) a direct sum of line bundles, or (b) isomorphic to $T_{\PP^2}(k)$ for some 
$k \in \Z$, where $T_{\PP^2}$ is the tangent bundle of $\PP^2$.
\label{uniform}
\end{thm}

Next let us introduce some definitions and results on line arrangements in $\PP^2$, to which we apply our main results. Let $\A$ be an arrangement of lines in 
$\PP^2$, namely, a finite set of lines in $\PP^2$. 
It can be naturally identified with a central arrangement $\overline{\A}$ of planes in $\C^3$.  Let 
$
L(\overline{\A}):=\{ \cap_{H \in \B} H \mid \B  \subset \overline{\A} \}
$
be the 
intersection lattice of $\overline{\A}$, with a partial order induced from the 
reverse inclusion, and let $\chi(\overline{\A};t)$ be the corresponding characteristic polynomial, see \cite{OT, DHA}. 

Then $\chi(\overline{\A};t)=\sum_{i=0}^3 (-1)^ib_{i}(\overline{\A}) t^{3-i}$,  where $b_i(\overline{\A})$ is the $i$-th Betti number of $M(\overline{\A})=V \setminus \cup_{H \in \overline{\A}} H$, see \cite{OT, DHA}. 
When $\A \neq \emptyset$, it is known that $\chi(\overline{\A};t)$ is divisible by 
$t-1$. Define $\chi(\A;t):=\chi(\overline{\A};t)/(t-1)$ and note that
$\chi(\A;t)=t^2-b_1(\A)t+b_2(\A)$, 
where 
$b_i(\A)$ is the $i$-th Betti number of $M(\A)=\PP^2 \setminus \cup_{H \in \A} H$, see \cite{OT, DHA}. Let us recall the definition of logarithmic vector fields and the freeness of 
arrangements.

\begin{definition}
Let $\alpha_H$ be a defining linear form for $H \in \A$. Then for 
$Q(\A):=\prod_{H \in \A} \alpha_H$, define 
$$
AR(\A):=AR(Q(\A)).
$$
For $H \in \A$, define 
\begin{eqnarray*}
AR_H(\A):&=&\{(a,b,c) \in S^{\oplus 3} 
\mid (a \partial_x+b \partial_y+c \partial_z)(\alpha_L) \in S\alpha_L\ (\forall L \in \A 
\setminus \{H\}),\ \\
&\ &(a \partial_x+b \partial_y+c \partial_z)(\alpha_H)=0\}.
\end{eqnarray*}
\end{definition}

The following is well-known, of which we give a proof for the 
completeness.

\begin{prop}
$AR(\A) \simeq AR_H(\A)$ for all $H \in \A$.
\label{basic}
\end{prop}

\proof
Let $\theta_E$ be the Euler derivation and 
$M:=AR(\A) \oplus S \theta_E$. 
Define a map $\varphi:
M \rightarrow AR_H(\A)$ by
$$
M \ni \theta \mapsto 
\theta - (\theta(\alpha_H)/\alpha_H)\theta_E \in
AR_H(\A).
$$
The kernel of $\varphi$ is clearly $S \cdot \theta_E$. Also, 
for any $\theta \in AR_H(\A)$, 
$\theta-(\theta(Q)/(\deg Q)Q)\theta_E 
\in AR(\A) \subset M$ is sent to $\theta$ by $\varphi$. Hence 
$$
AR_H(\A) \simeq M/S \cdot \theta_E \simeq AR(\A).
$$
\endproof

When $\A$ is free, i.e. when $AR(\A)$ is a free graded $S$-module, we have the following important result.

\begin{thm}[Terao's factorization, \cite{T2}]
Assume that $\A$ is free with $\exp(\A)=(d_1,d_2)$. Then 
$\chi(\A;t)=(t-d_1)(t-d_2)$.
\label{Teraofactorization}
\end{thm}

Note that Terao proved Theorem \ref{Teraofactorization}  in all dimensions, but the above case is enough for our purposes. Here is the nearly free version of this factorization result.

\begin{thm}[Factorization for nearly free arrangements, \cite{DSt15}]
Let $\A$ be nearly free with $\nexp(\A)=(d_1,d_2)$. Then 
$d_1+d_2=|\A|$ and 
$$
\chi(\A;t)=(t-d_1)(t-d_2+1)+1.
$$
\label{factorization}
\end{thm}

However, it is very difficult to determine whether a given arrangement is free 
or not in general, even for line arrangements. Here we recall a criterion for freeness. 
For that purpose, we need the following definition.

\begin{definition}
For a central arrangement $\A$ in $\C^3$ and $H \in \A$, define $\A^H:=
\{H \cap L \mid L \in \A \setminus \{H\}\}$ and 
$m^H(X):=|\{L \in \A \setminus \{H\} \mid 
L \cap H=X\}|$ for $X \in \A^H$. The pair $(\A^H,m^H)$ is 
called the Ziegler restriction of $\A$ onto $H$. Also, there is a canonical 
Ziegler restriction map
$$
\pi:AR_H(\A) \rightarrow D(\A^H,m^H),
$$
where 
$$
D(\A^H,m^H):=
\{(a,b) \in (S/\alpha_H)^{\oplus 2} \mid  
(a  \partial_x+b\partial_y)(\alpha_X) \in (S/\alpha_H) \alpha_X^{m^H(X)}\ 
(\forall X \in \A^H)\}
$$
and we choose $\alpha_H=z$ and hence $S/\alpha_H=\C[x,y]$.
\label{Zieglerrest}
\end{definition}

Note that the restriction $\A^H$ of a central arrangement $\A$ is also central. 
Since $D(\A^H,m^H)$ is also reflexive as an $S/\alpha_H$-module, it is free. If its free 
basis has degrees $(d_1,d_2)$, it is denoted as $\exp(\A^H,m^H)=(d_1,d_2)$. 
In this article, $d_1 \le d_2$ unless otherwise specified. Define 
$b_2(\A^H,m^H):=d_1d_2$. 

\begin{thm}[\cite{Z}]
Assume that $\A$ is free with $\exp(\A)=(d_1,d_2)$. Then 
$(\A^H,m^H)$ is also free with $\exp(\A^H,m^H)=(d_1,d_2)$, and the Ziegler 
restriction map $AR_H(\A) \rightarrow D(\A^H,m^H)$ 
is surjective for all $H \in \A$.
\label{ZE}
\end{thm}

The following is the arrangement version of Theorem \ref{Y2}. 

\begin{thm}[Yoshinaga's criterion, \cite{Y2}]
Let $\exp(\A^H,m^H)=(d_1,d_2)$. Then 
$b_2(\A) - d_1 d_2 \ge 0$, which coincides with 
$\dim \coker \pi$. Moreover, $\A$ is free with 
$\exp(\A)=(d_1,d_2)$ if and 
only if the equality holds.
\label{Y}
\end{thm}

In general, the pair $(\A,m)$, where $\A$ is a line arrangement, and 
$m:\A \rightarrow \Z_{>0}$ is called a multiarrangement. To investigate the exponents of the multiarrangement, the following easy 
lemma is important.

\begin{lem}[\cite{AN}, Lemma 4.2]
Let $(\A,m)$ be a multiarrangement in $\C^2$.
For $H \in \A$, let $\delta_H$ be a multiplicity such that 
$\delta_H(L)=1$ only when $H=L$, and $0$ otherwise. Then there is a homogeneos 
basis $\theta_1,\theta_2$ for $D(\A,m)$ such that 
$\alpha_H\theta_1,\theta_2$ form a basis for $D(\A,m+\delta_H)$.
\label{AN}
\end{lem}

Also, we use the following relation between Betti numbers and 
Chern classes, regarded as integers under the canonical identification
$H^{2i}(\PP^2,\Z)=\Z$, for details we refer to \cite[Proposition 5.18]{DS} or \cite[Corollary 4.3]{MuS}.
\begin{prop}
$b_i(\A)=(-1)^ ic_i(E_\A)$ for $i=0,1,2$.
\label{DS}
\end{prop}

\section{Proof of Theorem \ref{thm1}}
The following result is perhaps well-known, but we include a proof for reader's convenience.

\begin{prop} \label{p1} 
With the notation above, we have the following.
\begin{enumerate} 

\item  For any line $L$, one has $d_1^L +d_2^L=d-1$.

\item If the line $L_0$ is generic and the line $L$ arbitrary, then $d_1^{L_0} \geq d_1^{L}$.
In particular one has $d_1^{L_0} d_2^{L_0} \geq d_1^{L} d_2^{L}$.

\end{enumerate}  

\end{prop}

\proof

For the first claim, note that we have the following values for the first Chern classes, considered now as cohomology classes:
$$c_1(E_C)=c_1(T\langle C \rangle (-1))=-(d-1)\al,$$ see for instance \cite[ equation (3.2)]{DS14} and then obviously
$$c_1(E_C|_L) =c_1( \OO_L(-d_1^L) \oplus \OO_L(-d_2^L))=-(d_1^L +d_2^L)\be.$$
 Here $\al$ (resp. $\be$) are the canonical generators of $H^2(X,\Z)$ (resp. $H^2(L,\Z)$), with $X=\PP^2$. Let $i:L \to X$ denote the inclusion, and note that $i^*(\al)=\be$ and also $i^*(c_1(E_C))=c_1(i^*(E_C))$, with $i^*(E_C)=E_C|_L$.

For the second claim, see \cite[\S 1, Definition 2.2.3]{OSS}, but note that the ordering of the degrees $( d_1^L, d_2^L)$ in our paper is opposite from the ordering in \cite{OSS}. Indeed, the generic splitting type as defined in  \cite[\S 1, Definition 2.2.3]{OSS} corresponds to $d_2^L$ being minimal, hence in view of (1), to $d_1^L$ being maximal.

\endproof

The following result is the key step in proving Theorem \ref{thm1}

\begin{prop} \label{p2} Let  $L_0$ be a generic line in $\PP^2$.
With the above notation, we have the following.
\begin{enumerate} 

\item  For any line $L$ in $\PP^2$, one has 
$$\max(r-\dim N(f)_{r+d-3},0) \leq d_1^{L} \leq r,$$
 where $r=mdr(f)$.

\item If $d_1^{L_0} <(d-2)/2$, then $mdr(f) = d_1^{L_0}$.

\item If $d_1^{L_0}  \geq (d-2)/2$ and $d=2m$ is even, then $d_1^{L_0}=m-1$, $d_2^{L_0}=m$. In particular, this case  can occur for a free curve $C$ only if the exponents are $d_1=m-1$, $d_2=m$, while for a nearly free curve $C$ the exponents should necessarily be either $d_1=m-1$, $d_2=m+1$, or $d_1=d_2=m$.

\item If $d_1^{L_0}  \geq (d-2)/2$ and $d=2m+1$ is even, then $d_1^{L_0}=d_2^{L_0}=m$. In particular, this case can occur for a free curve  $C$ only if the exponents are $d_1=d_2=m$, while for a nearly free curve $C$ the exponents should necessarily be $d_1=m$, $d_2=m+1$. 

\end{enumerate}  

\end{prop}

\proof  To prove the inequality $d_1^{L} \leq r$ in claim (1), we use increasing induction on 
 $0 \leq k<d_1^L\leq d_2^L$, and prove  that $AR(f)_k=0$ in this range, using  the exact sequence \eqref{equa2}. Note that $AR(f)_0=0$ by our assumption $mdr(f) \geq 1$ in Introduction. To prove the  
 inequality $ r-\dim N(f)_{r+d-3}  \leq d_1^{L}$ in claim (1), we  use the exact sequence \eqref{equa2} for $k=r-1$, when we get
 $$0 \to H^0(L, \OO_L(r-1-d_1^{L})) \oplus H^0(L,\OO_L(r-1-d_2^{L})) \to N(f)_{r+d-3} \to \cdots.$$
If $d_1^{L}<r-\dim N(f)_{r+d-3}$, then $d_1^{L} \leq r-1$, and one must have
$$r-d_1^{L}=\dim H^0(L, \OO_L(r-1-d_1^{L})) \leq \dim N(f)_{r+d-3}.$$
This implies $d_1^{L} \geq  r-\dim N(f)_{r+d-3}$, and this is a contradiction. Finally note that $d_1^{L} \geq 0$, as follows from the last claim in Proposition \ref{key1}.

To prove the claim (2), take $k=d_1^{L_0}$, and note that the condition $d_1^{L_0} <(d-2)/2$ is equivalent to the condition $k+d-2 <T/2$, with $T=3(d-2)$. Now we use \cite[Corollary 4.3]{DPop} to conclude that the morphism $ N(f)_{k+d-2} \to  N(f)_{k+d-1} $ in the exact sequence \eqref{equa2}, which is induced by the multiplication by a generic linear form $\alpha_{L_0}$ defining ${L_0}$, is injective. Hence $AR(f)_k=H^0(L,\OO_L \oplus \OO_L(d_1^L-d_2^L)) \ne 0$,
which completes the proof for the second claim. 

The proof for the claims (3) and (4) goes along the same line, and we leave them to the reader.
\endproof

\noindent
\textit{Proof of Theorem \ref{thm1}}.
The claim (1) follows from the claim (1) in Proposition \ref{p2} (1), using \cite[Corollary 4.3]{DPop}.
The  inequality in (2) follows from Proposition \ref{p1}.
It is known that $c_2(E_C)=(d-1)^2-\tau(C)$, see  for instance \cite[equation (3.2)]{DS14}.
On the other hand, by  Theorem \ref{Y2} and the exact sequence \eqref{equa2} we get
$$c_2(E_C)-d_1^{L_0} d_2^{L_0}= \dim \ker \{ \alpha_{L_0} : N(f) \to N(f)\},$$
for a generic linear form $\alpha_{L_0}$. Using again \cite[Corollary 4.3]{DPop} it is clear that
$$\dim \ker \{\alpha_{L_0} : N(f) \to N(f)\}= \dim N(f)_{[T/2]}= \nu(C).$$
\owari
\medskip

By the proof above, we have the following:

\begin{cor}
$$c_2(E_C)-d_1^{L_0} d_2^{L_0}=(d-1)^2-\tau(C)-d_1^{L_0} d_2^{L_0}=\nu(C).
$$
\label{c2}
\end{cor}

\noindent
\textit{Proof of Corollary \ref{mainAG}}.
Combine Corollary \ref{c2} with Proposition \ref{dimNf} and the fact that 
the freeness is equivalent to $\nu(C)=0$. \owari
\medskip

Using Theorem \ref{thm1} (1) and Proposition \ref{p1} (1), we get the following result.
\begin{cor}
For any reduced plane curve $C:f=0$ of degree $d$, the image
$ \im (a_{E_C})$ is contained in
$$\{(r_0,d-1-r_0), (r_0-1,d-r_0),...,(r_0'+1, d-2-r_0'),(r_0', d-1-r_0')\},$$
where $r_0=\min (mdr(f), [(d-1)/2])$ and $r_0'=\max(mdr(f)-\nu(C),0)$. In particular, this set has at most 
$$r_0-r_0'+1 \leq \nu(C)+1$$
elements.
Moreover, if $C$ is nearly free with $\nexp(C)=(d_1,d_2),\ d_1 \le d_2$, and $L$ is any line, then $a_{E_{C}}(L)$ is either 
$(d_1-1,d_2)$ or $(d_1,d_2-1)$.
\label{keycor}
\end{cor}

The following result shows in particular that  $E_C$ is a uniform bundle for a nearly free curve $C$ if and only if the exponents of $C$ are equal.

\begin{cor}
With the above notation, assume that 
$c_2(E_C)=d_1(d_2-1)+1$, for some integers $1 \leq d_1 \leq d_2$. Then $C$ is nearly free with $\nexp(C)=(d_1,d_2)$ if and only if either
(1) $d_1=d_2=d'$ and
$a_{E_C}\equiv (d'-1,d')$, or (2) $d_1 <d_2$, 
$a_{E_C}$ is not constant and $\mbox{Im}(a_{E_C})=\{(d_1-1,d_2),(d_1,d_2-1)\}$.
\label{splittingtype}
\end{cor}

\proof
The ``if'' part follows from Corollary \ref{mainAG}. Assume that $C$ is 
nearly free with $\nexp(C)=(d_1,d_2)$. Then Corollary \ref{keycor} confirms that 
the splitting type is either $(d_1-1,d_2)$ or $(d_1,d_2-1)$. If $d_1=d_2$, then this is clearly
the case (1). Assume that $d_1 <d_2$, and that  
$a_E$ attains only the value $(d_1,d_2-1)$. Indeed, this value has to be in the image of $a_{E_C}$ due to Corollary \ref{mainAG}. Then Theorem \ref{uniform} 
says that, combining the fact that $C$ is nearly free, hence not free, 
$E \simeq T_{\PP^2}(k)$ for some $k \in \Z$. 
Then 
its splitting type is $(c-1,c)$, see \S 2.2, \cite{OSS}, thus $d_1+2=d_2$. 
Then the $S$-module $AR(f)$ is generated by 
one degree $d_1$-element and two degree $(d_1+2)$-elements. However, 
by the Euler sequence, the $S$-module of global sections 
of the twisted tangent bundle $T_{\PP^2}(k)$ is generated by three same degree elements, a contradiction.
\endproof

\section{Proof of Theorem \ref{mdr}}

First recall the definition of the local Tjurina number $\tau(C,p)$, where $p \in C$. Choose a local system of coordinates $(u,v)$ centered at $p$, and assume that the analytic germ $(C,p)$ is given by a local equation $g(u,v)=0$. Then one defines
$$\tau(C,p)= \dim_{\C} \frac{\OO_{\PP^2,p}}{(g,g_u,g_v)},$$
where $g_u,g_v$ are the partial derivatives of $g$ with respect to $u$ and $v$ respectively, and $(g,g_u,g_v)$ is the ideal spanned by these 3 germs in the local ring $\OO_{\PP^2,p}$ of analytic function germs at $p$.
Next recall 
the following basic result relating  $\tau(C)$ and $r=mdr(f)$ in  \cite[Theorem 3.2]{duPCTC}. 

\begin{prop}
\label{tjurina}
For any reduced plane curve of degree $d$, one has
$$\tau(d,r)_{min} \leq \tau(C) \leq \tau(d,r)_{max},$$
where $\tau(d,r)_{min}=(d-1)(d-r-1)$, and $\tau(d,r)_{max}=(d-1)(d-r-1)+r^2$ for $r \leq (d-1)/2$ and 
$$\tau(d,r)_{max}=(d-1)(d-r-1)+r^2- { 2r+2-d \choose 2},$$
for $(d-1)/2<r \leq d-1$.
\end{prop}
As already noted in \cite[Proof of Thm. 1.1]{D1}, the function 
$$r \mapsto \tau(d,r)_{max}$$
is strictly decreasing on the interval $[0,d-1]$.

Now we start the proof of Theorem \ref{mdr}. As we have seen in the previous section, one has $c_1=-(d-1)$ and $c_2=(d-1)^2-\tau(C)$.
Using this the equivalence of the claims (1) and (2) is clear, as both are equivalent to $\Delta=c_1^2-4c_2\geq 0$. Using the decreasing function
$\tau(d,r)_{max}$ and the remark that 
$$\tau(d,\frac{d-1}{2})_{max}=\frac{3}{4}(d-1)^2,$$
we see that (1) implies the inequality 
\begin{equation} \label{equa4.0} 
r \leq (d-1)/2. 
\end{equation}
Suppose now that (1) holds and replace $c_2$ by 
$$d_1^{L_0} d_2^{L_0}+\nu(C)=a(d-1-a)+\nu(C),$$
where we set $a=d_1^{L_0}$ for simplicity, and use Theorem \ref{thm1}.
The condition $\Delta=c_1^2-4c_2\geq 0$ now becomes
\begin{equation} \label{equa4.1} 
4a^2-4a(d-1)+(d-1)^2-4\nu(C) \geq 0.
\end{equation}
The associated equation
$$4a^2-4a(d-1)+(d-1)^2-4\nu(C) = 0,$$
has roots
$$a_1= \frac{d-1}{2}- \sqrt{ \nu(C)} \text{ and } a_2= \frac{d-1}{2}+\sqrt {\nu(C)}.$$
In view of the inequality \eqref{equa4.0}, it follows that the inequality \eqref{equa4.1} implies $a \leq a_1$. Hence we have shown that (1) implies (3). Conversely, if (3) holds, it follows that the inequality \eqref{equa4.1} holds, and hence $\Delta=c_1^2-4c_2\geq 0$.
Hence (3) implies (1) as well.

Now we show that (3) implies the equality $mdr(f)=d_1^{L_0}$.
If the curve $C$ is free, the claim is obvious. Otherwise, $\nu(C) \geq 1$, and hence the inequality (3) implies $d_1^{L_0} \leq (d-3)/2$. We conclude using Proposition \ref{p2} (2).

Finally, if $r=mdr(f) \leq (d-1)/4$, it follows that
$$\tau(C) \geq \tau(d,r)_{min} \geq \tau(d,\frac{d-1}{4})_{min}=\frac{3}{4}(d-1)^2.$$
In other words, the inequality $r=mdr(f) \leq (d-1)/4$ implies (2). 
To show that $mdr(f) \leq a_1$ implies (3), we use Proposition \ref{p2}, (1). This ends the proof of Theorem \ref{mdr}.

Combining Corollary \ref{keycor} with Theorem \ref{mdr}, we have the following.

\begin{cor}
 If the Chern polynomial $c_{E_C}(t)=1+c_1t+c_2 t^2\in \Z[t]$ of $E_C$ has real roots, then
$$
s_C-1\le \nu(C) \le \left(
\frac{d-1}{2}-mdr(f)\right)^2,
$$
where $s_C=|\im (a_{E_C})|$.
\label{ineq}
\end{cor}

\begin{ex} \label{ZieglerEx}
Let $p_1, \dots, p_6\in\PP^2$ be $6$ points such that no three are colinear. 
Let $\A=\{H_1, \dots, H_9\}$ be the edges of the corresponding hexagon and three diagonals (more precisely, the lines $p_1p_2, p_2p_3, p_3p_4, p_4p_5, p_5p_6, p_6p_1, p_1p_4, p_2p_5$ and $p_3p_6$), such that each line $H_j$ contains exactly 2 triple points and  4 nodes. 

Denote by $\A:f=0$ (resp. by $\A':f'=0$) the corresponding line arrangement when the 6 vertices of the hexagon are (resp. are not) on a conic. 
Then it is known that this pair of Ziegler's line arrangements has the following properties, see for instance \cite[Remark 8.5]{DHA}.

\begin{enumerate} 

\item  Both arrangements have the same intersection lattice. In particular, they have 18 double points and 6 triple points. It follows that $\tau(\A)=\tau(\A')=42<(3/4)\cdot 64=48.$ Hence the Chern polynomials 
$$c_{E_{\A}}(t)=c_{E_{\A'}}(t)=1-8t+22t^2$$
do not have real roots in view of Theorem \ref{mdr}.

\item $r=mdr(f)=5$ and $r'=mdr(f')=6.$

\item Since $d_1^{L_0} \leq (d-1)/2=4$, it follows that we are in the case (4) of Proposition \ref{p2}. Hence $d_1^{L_0}=d_2^{L_0}=4$ for both arrangements $\A$ and $\A'$. In particular $d_1^{L_0} \ne mdr(f)$ in these two cases.

\end{enumerate} 
Note that one can consider a family $\A'_t:f_t=0$ of arrangements as above, where the 6 vertices of the hexagon are not on a conic, for $t\ne0$, degenerating at an arrangement $\A'_0=\A$, where the 6 vertices of the hexagon are on a conic. Note that under this degeneration the invariant $mdr(f_t)$ drops by one, while the corresponding splitting invariant $d_{1,t}^{L_0}$ stay constant.

\end{ex}

\section{Application to line arrangements}

In this section let us apply the main results in this paper to 
the case when $C$ is a finite set of lines in $\PP^2$, i.e., 
an arrangement of lines. 
For that purpose, let us show the following generalization 
of Ziegler's result in \cite{Z}. Only in this result $\ell$ is arbitrary. 

\begin{thm}
Let $\K$ be an arbitrary field and 
$\A$ be a central arrangement in $V=\K^\ell$ and $S:=\K[x_1,\ldots,x_\ell]$. 
Let $\pi:AR_H(\A) \rightarrow D(\A^H,m^H)$ be the Ziegler restriction map. 
Let $\theta_1,\ldots,\theta_s \in AR_H(\A)$ satisfy that 
$\pi(\theta_1),\ldots,\pi(\theta_s)$ generate $\mbox{Im}(\pi)$ as an $(S/\alpha_H)$-module. Then 
$\theta_1,\ldots,\theta_s$ generate $AR_H(\A)$ as an $S$-module.
\label{Ziegler}
\end{thm}

\proof
Let the images of $\theta_1,\ldots,\theta_s$ by $\pi$ generate $M:=\mbox{Im} (\pi)$, $\alpha_H=x_\ell$ and let $0 \neq \theta \in AR_H(\A)$ be a minimal degree element. 
Then we may show that $\pi(\theta) \neq 0$. Assume not. Then $\theta=x_\ell \theta'$ for 
some $\theta' \in \Der S$. Since 
$\varphi(\alpha_H)=0$ for $\varphi \in AR_H(\A)$, it follows that 
$0 \neq \theta' \in AR_H(\A)$ with $\deg \theta'< \deg \theta$, a contradiction. 
Hence $0 \neq \pi(\theta)=\sum_{i=1}^s a_i \pi(\theta_i)$ 
for $a_i \in \K$ such that $a_i=0$ if 
$\deg \theta_i> \deg \theta$. Then $\theta-\sum_{i=1}^s a_i \theta_i=x_\ell \theta'$. By the 
same reason, $\theta' \in AR_H(\A)$ whose degree is strictly lower than that of 
$\theta$, a contradiction. Hence $\theta'=0$, hence the lowest degree derivations in 
$AR_H(\A)$ can be expressed by $\theta_1,\ldots,\theta_s$.  

Now assume that the statement holds true for homogeneous derivations in $AR_H(\A)$ whose 
degree is less than $d$. Since $AR_H(\A)$ is graded, it suffices to show the statemet for 
homogeneous parts. Let $\theta \in AR_H(\A)_k$. If $x_\ell \mid \theta$, then 
apply the induction hypothesis to $\theta/x_\ell \in AR_H(\A)_{k-1}$, which 
completes the proof. 
Assume not. Then the same argument as above implies 
that $\theta-\sum_{i=1}^s f_i \theta_i =x_\ell \theta'$ for some $\theta' \in AR_H(\A)_{k-1}$.
Again the induction hypothesis implies that $\theta'=\sum_{i=1}^s g_i \theta_i$, which 
completes the proof. 
\endproof

\begin{rk}
When $\pi$ is surjective, $M=D(\A^H,m^H)$ in terms of Theorem \ref{Ziegler}. Hence the classical result by Ziegler in \cite{Z} asserting that 
$\pi$ is surjective if $\A$ is free can be regarded as 
a special case of Theorem \ref{Ziegler}.
\end{rk}

In general, it is very difficult to investigate 
the splitting type of vector bundles onto projective lines. Contrary to it, 
for line arrangements, we can use the technique of multiarrangements to do it.
What makes this analysis work well is the following exact sequence which 
hold true when $\ell=3$:
$$
0 \rightarrow \widetilde{AR_H(\A)} \stackrel{\cdot \alpha_H}{\rightarrow}
\widetilde{AR_H(\A)} \rightarrow \widetilde{D(\A^H,m^H)} \rightarrow 0.
$$
This follows by Proposition \ref{key1} and the fact that $\ell=3$. This implies the 
isomorphism:
$$
\widetilde{D(\A^H,m^H)} \simeq \widetilde{AR_H(\A)}|_H.
$$
Hence to know $a_E(H)$ for $H \in \A$ is the same as to know 
$\exp(\A^H,m^H)$. We use this isomorphism frequently in the rest of this article.

From now on, we assume that $\ell=3$, so every arragement is that of lines in 
$\PP^2$.
Next we prove 
the nearly 
free version of Corollary \ref{keycor}.

\begin{thm}
Assume that the line arrangement $\A$ is neary free with $\nexp(\A)=(d_1,d_2)$. Then 
$\exp(\A^H,m^H)=(d_1-1,d_2)$ or $(d_1,d_2-1)$.
\label{main2}
\end{thm}

\proof
Apply Corollary \ref{keycor} and the exact sequence above. 
\endproof

\begin{ex}
Let $\A=\{xyz(y-z)(x+2y+3z)=0\}$. Then 
$\chi(\A;t)=t^2-4t+5=(t-2)^2+1$. It is easy to check that (e.g., 
use Theorem \ref{criterion} below) $\A$ is nearly free with $\nexp(\A)=(2,3)$. Also, 
$\exp(\A^H,m^H)=(2,2)$ and 
$\exp(\A^L,m^L)=(1,3)$ for $H:z=0$ and $L:x+2y+3z=0$. Hence both case in Theorem \ref{main2}  
can occur in general. 
\end{ex}

Similarly to Theorem \ref{Y}, 
we may give a sufficient condition for a line arrangement to be nearly free 
following Theorem \ref{mainAG}.

\begin{thm}[Near freeness condition]
Let $\chi(\A;t)=(t-d_1)(t-d_2+1)+1$ with $d_1 \le d_2$. Then 
$\A$ is nearly free if there is $H \in \A$ such that 
$\chi(\A;0)-b_2(\A^H,m^H)=1$. Here, $b_2(\A^H,m^H)$ is the product of 
exponents of $(\A^H,m^H)$, see \cite{ATW}. In particular, 
$\nexp(\A)=(d_1,d_2)$ or $(d_1+1,d_2-1)$ if $d_1=d_2$ or $d_1+2=d_2$, and 
$\nexp(\A)=(d_1,d_2)$ otherwise. 
\label{criterion}
\end{thm}

\proof
Apply Theorem \ref{mainAG}, 
Proposition \ref{p2} and the isomorphism above. \endproof
\medskip

\begin{rk}
We do not know whether the result like Corollary \ref{splittingtype} holds true for 
$H \in \A$. In other words, we do not know whether for a nearly free arrangements, 
there is $H \in \A$ such that $\exp(\A^H,m^H)=(d_1,d_2-1)$. We do not have 
any counter example to this statement. 
\end{rk}

Next let us study the addition-deletion type results for free and nearly free 
arrangements.

\begin{thm}
(1)\,\, 
Let $\A$ be free with $\exp(\A)=(d_1,d_2)$ with $d_1 \le d_2$. Let $L \not \in \A$ be a line. Then 
$\B:=\A \cup\{L\}$ is nearly free if 
$|\B^L|=d_2+2$.

(2)\,\, 
Let $\A$ be free with $\exp(\A)=(d_1,d_2)$ with $d_1 \le d_2$. Let $H \in \A$ be a line. Then 
$\B:=\A \setminus \{H\}$ is nearly free if 
$|\A^H|=d_1$.
\label{addition2}
\end{thm}

\proof
(1)\,\, 
Assume that $d_1=d_2=:d$. Then $\exp(\B^H,m^H)=(d,d+1)$ for any $H \in \B$ by the 
argument in \cite{A}. Also, $\chi(\B;0)=d^2+d+1$. Hence their difference is $1$, so Theorem \ref{criterion} completes the proof. 

So we may assume that 
$d_1<d_2$. By Theorem \ref{ZE} and Lemma \ref{AN}, 
$\exp(\B^H,m^H)$ are either $(d_1+1,d_2)$ or $(d_1,d_2+1)$ for $H \in \A$. 

Case 1. Assume that 
there is $H \in \B$ such that 
$\exp(\B^H,m^H)=(d_1+1,d_2)$.  Then for $\pi:AR_H(\B) \rightarrow D(\B^H,m^H)$, 
$\dim \coker \pi =d_1d_2+d_2+1-(d_1+1)d_2=1$. Hence Theorem \ref{criterion} 
completes the proof. 

Case 2. Assume that $\exp(\B^H,m^H)=(d_1,d_2+1)$ for all $H \in \B$. 
We may assume that $H \neq L$. 
Since 
$b_2(\B)=d_1d_2+d_2+1$, it holds that 
$\dim \coker \pi =d_2-d_1+1$, hence $\B$ is not free by Theorem \ref{Y}. 
Let $\theta_1,\theta_2$ be a basis for 
$AR_H(\A)$ with $\deg \theta_i=d_i$. Then clearly 
$\alpha_L \theta_i \in AR_H(\B)_{d_i+1}$, and $\alpha_L\theta_1,\alpha_L \theta_2$ are 
$S$-independent. Let $\eta_1,\eta_2$ be a basis for $D(\B^H,m^H)$ with 
$\deg \eta_1=d_1,\deg \eta_2=d_2+1$. If $AR_H(\B)_{d_1} \neq 0$, then clearly 
$AR_H(\B)$ is free, which is a contradiction. Hence 
$\eta_1 \not \in \mbox{Im} \pi$. So we may assume that, by putting $\alpha_H=z$, 
$\pi(\alpha_L \theta_1)=x \eta_1$. 
Since $\theta_1$ and $\theta_2$ form a basis for $AR_H(\A)$, $\pi(\theta_1)$ and 
$\pi(\theta_2)$ are $(S/\alpha_H)$-independent. 
So we may assume that 
$\pi(\alpha_L \theta_2)=\eta_2$.  
Now taking $\dim \coker \pi$ into account, there has to be  
$\theta_3 \in AR_H(\B)_{d_2+1}$ such that $\pi(\theta_3)=y^{d_2-d_1+1}\eta_1$.
Then by Theorem \ref{Ziegler}, $\alpha_L\theta_1,\alpha_L\theta_2,\theta_3$ satisfy the condition for 
$\B$ to be nearly free.

(2) Assume that $\A$ is free with $\exp(\A)=(d_1,d_2)_\le$. By the deletion-restriction, 
$\chi(\B;t)=(t-d_1)(t-d_2+1)+1$. Let $L \in \B$. Since 
$\exp(\A^L,m^L)=(d_1,d_2)$ by Theorem \ref{ZE}, $\exp(\B^L,m^L)$ is either 
$(d_1-1,d_2)$ or $(d_1,d_2-1)$. If the latter, then Theorem \ref{criterion} 
completes the proof. Assume the former for all $L \in \B$. Then 
$$
b_2(\B)-(d_1-1)d_2=d_2-d_1+1.
$$
By Lemma \ref{AN}, there are a basis $\theta_1,\theta_2$ for $D(\B^L,m^L)$ such that 
$\deg \theta_1=d_1-1,\ \deg \theta_2=d_2$ and 
$x \theta_1,\theta_2$ form a basis for $D(\A^L,m^L)$, where $\alpha_H:=x$. Since 
$\pi:AR_L(\A) \rightarrow D(\A^L,m^L)$ is surjective by 
Theorem \ref{ZE}, there are derivations 
$\varphi_1,\varphi_2 \in AR_L(\A)$ such that $\pi(\varphi_1)=x\theta_1$ and 
$\pi(\varphi_2)=\theta_2$. Since $AR_L(\A) \subset AR_L(\B)$ and 
$\dim \coker (\pi_\B:AR_L(\B) \rightarrow D(\B^L,m^L))=d_2-d_1+1$, there 
is a derivation $\varphi_3 \in AR_L(\B)_{d_2}$ such that 
$\pi(\varphi_3)=y^{d_2-d_1+1}\theta_1$. Hence 
$$
\mbox{Im}\pi_\B =\langle x \theta_1,y^{d_2-d_1+1} \theta_1,\theta_2\rangle_{S/\alpha_L}.
$$
Since $\deg \varphi_3=d_2$ and clearly there is a relation among 
$\varphi_1,\varphi_2,\varphi_3$ at degree $d_2+1$, Theorem \ref{Ziegler} implies that 
$\B$ is nearly free.\endproof

The following is a nearly free version of the results in \cite{A}.

\begin{thm}
Let $\A$ be an arrangement of lines in $\PP^2$ with 
$\chi(\A;t)=t^2-b_1t+b_2$, where $b_1=|\A|-1$. Let 
$\chi(\A;t)=(t-a)(t-b)+1$ with real number $a\le b,\ a+b=b_1$. Then 
$\A$ is nearly free if there is $H \in \A$ such that 
\begin{itemize}
\item[(1)]
$|\A^H|=b+1$, or 
\item[(2)]
$|\A^H|=a+1$ and $b\neq a+2$.
\end{itemize} 
\label{sing}
\end{thm}

\proof Immediate from Theorem \ref{criterion} and the argument in \cite{A}.
\endproof

Now let us apply the results in this paper to 
show near freeness of some line arrangements.

\begin{ex}
(1)\,\,
Let $\A$ be defined by 
$$
xz(x^2-y^2)(x^2-2y^2)(y-z)=0..
$$
Then it is easy to check (see \cite{OT} for example) to show that 
$\chi(\A;t)=(t-3)^2=(t-2)(t-4)+1$, but $\A$ is not free. We can check the non-freeness 
by several way, here we use Theorem \ref{criterion}. It is easy to check that 
$\exp(\A^H,m^H)=(2,4)$ for any $H$ going through the origin. Hence 
Theorem \ref{criterion} implies that $\A$ is nearly free with $\nexp(\A)=(2,5)$.

(2)\,\,
Let $\B$ be defined by 
$$
xyz(x^2-z^2)(y^2-z^2)(x-y+z)(x-y-z)=0.
$$
Then 
$\chi(\B;t)=(t-4)^2+1$. 
Also, 
$|\B^H|=5=4+1$ for $\ker \alpha_H=x-y\pm z$. Hence 
Theorem \ref{sing} implies that $\B$ is nearly free with 
$\nexp(\B)=(4,5)$.
\end{ex}

\begin{thm}[Addition theorem for free and nearly free arrangements]
Let $\A$ be an arrangement in $\PP^2$, $H \in \A$ and 
let 
$\B:=\A \setminus \{H\}$. Also, let $d_1 \le d_2$ be two 
non-negative integers. Then the two of the following three 
implies the third:

(1)\,\,
$\A$ is nearly free with $\nexp(\A)=(d_1+1,d_2+1)$.

(2)\,\,
$\B$ is free with $\exp(\B)=(d_1,d_2)$.

(3)\,\,
$|\A^H|=d_2+2$.
\label{addition3}
\end{thm}

\proof
(1) and (2) implies (3) by 
two factorizations Theorems \ref{Teraofactorization} and \ref{factorization} and the 
deletion-restriction formula. 
If we assume 
(2) and (3), then Theorem \ref{addition2} (1) implies (1). 

Assume that (1) and (3). 
Let $L \in \B$. Then Theorem \ref{main2} shows that 
$\exp(\A^L,m^L)=:\exp(\A'',m)=(d_1,d_2+1)$ or 
$(d_1+1,d_2)$. Hence $\exp(\B^L,m^L)
=:\exp(\B'',k)$ could be one of $(d_1,d_2),\ 
(d_1+1,d_2-1)$ or $(d_1-1,d_2+1)$. Note that 
$$
\chi(\B;0)=b_2(\B)=d_1d_2 \ge b_2(\B'',k)
$$ 
by Theorem \ref{Y} and the deletion-restriction. Hence the proof is completed if 
$\exp(\B'',k)=(d_1,d_2)$. Assume that $\exp(\B'',k)=(d_1+1,d_2-1)$. Then by Theorem 
\ref{Y}, 
$$
b_2(\B)-b_2(\B'',k)=-d_2+d_1+1 \ge 0. 
$$
Hence $d_2=d_1+1$ or $d_2=d_1$. For the former, 
Theorem \ref{Y} confirms that $\B$ is free with 
exponents $(d_1,d_1+1)$.  For the latter, 
$b_2(\B)-b_2(\B'',k)=1$. Since 
$b_2(\B)=d_1^2=(d_1+1)(d_1-1)+1$, Theorem \ref{criterion} shows that 
$\B$ is nearly free with $\nexp(\B)=(d_1-1,d_1+2)$. Hence 
$0 \neq \alpha_H AR_L(\B)_{d_1-1}\subset AR_L(\A)_{d_1}=(0)$,  
a contradiction. 

Now assume that $\exp(\B'',k)=(d_1-1,d_2+1)$. Then $b_2(\B)-
b_2(\B'',k)=d_2-d_1+1$. This occurs only when 
$\exp(\A^L,m^L)=(d_1,d_2+1)$. Let $\alpha_H=y$. Then by Lemma \ref{AN}, 
we may choose a basis 
$\theta_1,\theta_2$ for $D(\B'',k)$ with $\deg \theta_1=d_1-1,\ 
\deg \theta_2=d_2+1$ such that $y \theta_1,\theta_2$ form a basis for 
$D(\A^L,m^L)$. Note that 
$b_2(\A)-b_2(\A'',m)=d_2-d_1+1$ too. Since $\nexp(\A)=(d_1+1,d_2+1)$, 
there is $\varphi \in AR_L(\A)_{d_2+1}$ such that $\pi(\varphi)=\theta_2$, and 
there are $\psi_1,\psi_2 \in AR_L(\A)$ such that 
$\pi(\psi_1)=xy \theta_1$ and $\pi(\psi_2)=y^{d_2-d_1+2}\theta_1$ for 
the Ziegler restriction map $\pi:AR_L(\A) 
\rightarrow D(\A^L,m^L)$. 
Hence $\coker \pi$ has a basis 
$$
y\theta_1,y^2\theta_1,\ldots,y^{d_2-d_1+1} \theta_1,
$$
whose dimension is surely $d_2-d_1+1$. 
Now consider the basis of $\coker (\pi':AR_L(\B) \rightarrow D(\B'',k))$. 
Assume that $\coker \pi' \ni y^i \theta_1$, where this $i$ is the largest one 
satisfying this. Then 
\begin{eqnarray*}
y^{i+1} \theta_1 
\in \mbox{Im}\pi'\ \mbox{and}\ 
y^{i} \theta_1 
\not \in \mbox{Im} \pi '&\iff&
y^{i+2} \theta_1 
\in \mbox{Im}\pi\ \mbox{and}\ 
y^{i+1} \theta_1 
\not \in \mbox{Im} \pi \\
&\iff&
0 \neq y^{i+1} \theta_1 
\in \mbox{coker}\pi'.
\end{eqnarray*}
Hence $d_2-d_1=i$, and 
$\coker \pi'$, whose dimension is $d_2-d_1+1$, could contain at most 
$$
\theta_1,\ x \theta_1,\ldots,x^{d_2-d_1+1} \theta_1,\ 
y\theta_1, y^2\theta_1,\ldots,y^{d_2-d_1} \theta_1.
$$
Hence $x \theta_1 \in D(\B'',k)$ since $\B$ is not free for 
$b_2(\B)=d_1d_2 \neq (d_1-1)(d_2+1)$. 

Summarizing, $AR_L(\B)$ has a generator 
$\eta_1$ of degree $d_1$, $\eta_2$ of degree $d_2$ and 
$\xi$ of degree $d_2+1$ such that $\pi'(\eta_1)=x \theta_1,\ 
\pi'(\eta_2)=y^{d_2-d_1+1} \theta_1,\ \pi'(\xi)=\theta_2$. We show that this cannot occur.
By the choice of $\eta_1,\eta_2$, 
there 
is $\eta \in AR_L(\B)_{d_2}$ such that 
$$
y^{d_2-d_1+1} \eta_1-x \eta_2=z \eta,
$$
where we set $z=\alpha_L$.
Since $AR_L(\B)_{ \le d_2}$ is generated by $\{\eta_1,\ \eta_2\}$, 
there are $a \in \K$ and $ g\in S_{d_2-d_1}$ such that 
$$
\eta= g\eta_1+a\eta_2 \iff 
(y^{d_2-d_1+1}-zg)\eta_1=(x+az )\eta_2.
$$
Since $y^{d_2-d_1+1}-zg$ and $x+az$ are coprime, there is $\eta_0 \in \Der(S)_{d_1-1}$ 
such that 
$$
\eta_0=\displaystyle \frac{\eta_1}{x+az}=\displaystyle \frac{\eta_2}{y^{d_2-d_1+1}-zg}.
$$
By the former expression of $\eta_0$, to show 
$\eta_0 \in AR_L(\B)_{d_1-1}$, it suffices to show that $\eta_0(x+az)\in 
S \cdot (x+az)$ if $\ker (x+az) \in \B$. 
Since $(x+az) \nmid (y^{d_2-d_1+1}-zg)$ and $\eta_2(x+az) \in S \cdot (x+az)$, 
it follows that 
$$ 
\eta_0(x+az)=\displaystyle \frac{\eta_2(x+az)}{y^{d_2-d_1+1}-zg}
\in S (x+az),
$$
which implies $0 \neq \eta_0 \in AR_L(\B)_{d_1-1}=(0)$, a contradiction. \endproof

\begin{thm}[Deletion theorem for free and nearly free arrangements]
Let $\A$ be an arrangement in $\PP^2$, $H \in \A$ and 
let 
$\B:=\A \setminus \{H\}$. Also, let $d_1 \le d_2$ be two 
non-negative integers. Then the two of the following three 
implies the third:

(1)\,\,
$\A$ is free with $\exp(\A)=(d_1,d_2)$.

(2)\,\,
$\B$ is nearly free with $\nexp(\B)=(d_1,d_2)$.

(3)\,\,
$|\A^H|=d_1$.
\label{addition4}
\end{thm}

\proof As in the proof of Theorem \ref{addition3}, it suffices to show that 
(2) and (3) imply (1). 
By Theorem \ref{Teraofactorization} and 
the deletion-restriction theorem, $\chi(\A;t)=(t-d_1)(t-d_2)$ and 
$\chi(\B;t)=(t-d_1)(t-d_2+1)+1$. Since 
$(a,b)_\le:=\exp(\B^L,m^L)$ is either $(d_1-1,d_2)$ or $(d_1,d_2-1)$ by Theorem \ref{main2} for $L \in \B$, 
$\exp(\A^L,m^L)$ is either $(d_1,d_2),\ (d_1-1,d_2+1)$ or $(d_1+1,d_2-1)$ by 
Lemma \ref{AN}.  
If it is $(d_1,d_2)$, then Theorem \ref{Y} completes the proof. Assume that 
$\exp(\A^L,m^L)=(d_1+1,d_2-1)$, then 
$$
\dim \coker (\pi_L:AR_L(\A) 
\rightarrow D(\A^L,m^L))=-d_2+d_1+1 \ge 0.
$$
Since $d_1 \le d_2$, we have $d_1=d_2$ or $d_1+1=d_2$. Assume that 
$d_1=d_2$. Then Theorem \ref{criterion} shows that $\A$ is nearly free with 
$\nexp(\A)=(d_1-1,d_1+2)$. Hence 
$(0) \neq AR_L(\A)_{d_1-1} \subset AR_L(\B)_{d_1-1}=(0)$, a contradiction. 
If 
$d_2=d_1+1$, then Theorem \ref{Y} says that $\A$ is free.

Hence we may assume that $\exp(\A^L,m^L)=(d_1-1,d_2+1)$ for all $L \in \B$. If $d_1=d_2$, then the above completes the proof. So we may assume that  
$d_1 < d_2$. This occurs only when 
$\exp(\B^L,m^L)=(d_1-1,d_2)$ for all $L \in \B$. Hence putting $\alpha_H:=x$, 
there is a basis 
$\theta_1,\theta_2$ of degree $d_1-1,d_2$ for $D(\B^L,m^L)$ such that 
$ \theta_1,x\theta_2$ form a basis for $D(\A^L,m^L)$. 
Note that, for $\pi_\B:AR_L(\B) \rightarrow D(\B^L,m^L)$, 
$$
\dim \coker_\B \pi=d_2-d_1+1.
$$
This also holds true for $\pi_\A:AR_L(\A) \rightarrow D(\A^H,m^H)$ by the assumption 
of the 
exponents. 
Let $\varphi_1,\varphi_2, \varphi_3$ be a generator for the 
nearly free module $AR_L(\B)$ such that 
$\deg \varphi_1=d_1,\deg \varphi_2=\deg \varphi_3=d_2$ and 
$\pi_\B(\varphi_1)=\alpha \theta_1,\ 
\pi_\B(\varphi_2)=\beta^{d_2-d_1+1}\theta_1$ and 
$\pi_\B(\varphi_3)=\theta_2$ for $\C$-independent linear forms 
$\alpha,\beta,z$. 

Since $x \varphi_3 \in AR_L(\A)$ and $ \pi_\A(x\varphi_3)=x \theta_2$, $\dim \coker \pi_\A$
comes from $\M_\A:=\mbox{Im} \pi_\A/\theta_1 \subset M_\B$, and they have the same codimension. Hence $\theta_1,\theta_2 \in AR_L(\A)$. 
 
Hence $\varphi_1,\varphi_2,x \varphi_3$ generates $AR_L(\A)$ by Theorem \ref{Ziegler}. 
Since $\pi_\A(\beta^{d_2-d_1+1}\varphi_1-\alpha \varphi_2)=0$, the same argument 
as the above shows that there is $\varphi \in AR_L(\A)_{<d_1} \subset 
AR_L(\B)_{<d_1}$, a contradiction. \endproof

\section{Proof of Theorem \ref{456}}

Let us prove Theorem \ref{456}.
If $|\A^H| \le 2$, there is nothing to show. Assume that 
$|\A^H|=3$. Then by \cite{W} with the assumption on $H$, the exponents of the Ziegler restriction 
onto $H$ is combinatorial, and it is a generic splitting type. Hence 
Corollary \ref{mainAG} completes the proof. 

Next assume that $|\A^H|=4$. 
Then by Theorem 1.6 in \cite{A1} with the assumption on $H$, $\exp(\A^H,m^H)=(e_1,e_2)$ is either 
$(d,d),\ (d,d+1)$ or $(d,d+2)$. For the former two cases, they are generic 
splitting types. Hence Corollary \ref{mainAG} and Theorem \ref{Y} says that, 
if $\chi(\A;t)$ has integer roots, then they are either free or nearly free 
if the roots satisfy conditions in Corollary \ref{mainAG} or Theorem \ref{Y}, and 
not either otherwise. So the rest case is when $|\A|=1+2d,\ 
(e_1,e_2)=(d-1,d+1)$ and $b_2(\A)=d^2$. In this case, Corollary \ref{mainAG} shows that 
$\A$ is nearly free. 
\medskip

\begin{rk}
If $H$ in Theorem \ref{456} contains a 
high multiplicity point, then freeness depends only on $L(\A)$, see 
Proposition 7.4 in \cite{WY} or Corollary 1.4 in \cite{D1}  for example. 
Also in this case, the following proposition shows the same holds true for 
near freeness.
\end{rk}

\begin{prop}
Let $\A$ be a line arrangement such that 
$(\A^H,m^H)$ is not balanced for some $H \in \A$, i.e., 
there is $X \in \A^H$ such that 
$2m^H(X) \ge |m^H|:=\sum_{Y \in \A^H} m^H(Y)=|\A|-1$. Then the near freeness of $\A$ depends only 
on $L(\A)$.
\label{notbalanced}
\end{prop}

\proof
Since $(\A^H,m^H)$ is not balanced, 
$\exp(\A^H,m^H)=(e_1,e_2)_\le$ with 
$m^H(X)=e_2$ (See \cite{Y3}, p 11 for instance). If $b_2(\A)=e_1e_2$, then $\A$ is free by 
Theorem \ref{Y}, hence not nearly free. So assume not, then $\A$ is not free.  
If $b_2(\A)-e_1 e_2=1$, then $\A$ is nearly free by Theorem \ref{criterion}. 
By Theorem \ref{factorization}, for $\A$ to be nearly free, 
there exist $d_1 \le d_2$ such that 
$$
\chi(\A;t)=(t-d_1)(t-d_2+1)+1.
$$
Also, by Corollary \ref{splittingtype}, 
if $(e_1,e_2)$ is neither $(d_1-1,d_2)$ nor $(d_1,d_2-1)$, then 
$\A$ is not nearly free. If $(e_1,e_2)=(d_1,d_2-1)$, then $\A$ is nearly free as shown 
above. Now assume that $(e_1,e_2)=(d_1-1,d_2)=(|\A|-1-m^H(X),m^H(X))$. Also, 
we may assume that $d_1 < d_2$. 
Let $H=\ker z$ and $X=H \cap \{y=0\}$. Then 
$\theta_1:=(\prod_{X \not \subset L \in \A} \alpha_L) \partial_x \in AR(\A)_{d_1-1}.
$ Hence by the Ziegler restriction map $\pi_H:AR_H(\A) \rightarrow D(\A^H,m^H)$, 
$\theta_1$ goes to $\varphi_1$, where $\varphi_1,\varphi_2$ form a basis for 
$D(\A^H,m^H)$ of degree $e_1,e_2$ respectively. Since 
$d_1 < d_2$, $(e_1,e_2)=a_{E_\A}(H)=
(d_1-1,d_2)$ is a generic splitting type of $E_\A$, which does not coincide with 
$(d_1,d_2-1)$. Hence Corollary \ref{mainAG} shows that $\A$ is not 
nearly free.
\endproof


\begin{thebibliography}{00}

\bibitem{A1}
T. Abe, 
Chambers of 2-affine arrangements and freeness of 3-arrangements. 
\textit{J. Alg. Combin}., \textbf{38} (2013), no. 1, 65--78. 

\bibitem{A}
T. Abe, 
Roots of characteristic polynomials and 
and intersection points of line arrangements. 
\textit{J. Singularities}, 
\textbf{8} (2014), 100--117.

\bibitem{AN}
T. Abe and Y. Numata, 
Exponents of $2$-multiarrangements and multiplicity lattices. 
\textit{J. Alg. Combin.}, \textbf{35} (2012),  
no. 1, 1--17.


\bibitem{ATW}
T. Abe, H. Terao and M. Wakefield, 
The characteristic polynomial of a multiarrangement.
\textit{Adv. in Math.}, 
\textbf{215} (2007), 825--838.



\bibitem{DS}
G. Denham and M. Schulze,
Complexes, duality and Chern classes of logarithmic forms along hyperplane arrangements.
\textit{Arrangements of hyperplanes—Sapporo} 
2009,  27-–57, Adv. Stud. Pure Math., 62, Math. Soc. Japan, Tokyo, 2012.

\bibitem{D1}
A. Dimca, Curve arrangements, pencils, and Jacobian syzygies,  \textit{Michigan Math. J.} \textbf{66} (2017), 347--365.

\bibitem{D2}
A. Dimca, 
Freeness versus maximal global Tjurina number for plane curves.
\textit{Math. Proc. Cambridge Phil. Soc.}, \textbf{163} (2017), 161--172.


\bibitem{DHA}  A. Dimca,   {\em Hyperplane Arrangements: An Introduction}, Universitext, Springer, 2017.

\bibitem{DPop} A. Dimca, D. Popescu, 
Hilbert series and Lefschetz properties of dimension one almost complete intersections, \textit{Comm. Algebra} \textbf{44} (2016), 4467--4482.


\bibitem{DS14} A. Dimca, E. Sernesi,  Syzygies and logarithmic vector fields along plane curves,
\textit{Journal de l'\'Ecole polytechnique-Math\'ematiques} \textbf{1} (2014), 247-267.


\bibitem{DSt15}
A. Dimca and G. Sticlaru, 
Nearly free divisiors and rational cuspidal curves.
arXiv:1505.00666.


\bibitem{duPCTC} A.A. du Plessis,  C.T.C. Wall, Application of the theory of the
discriminant to highly singular plane curves, \textit{Math. Proc. Camb.
Phil. Soc.},  \textbf{126} (1999), 259-266. 






\bibitem{FV1}  D. Faenzi, J. Vall\`es, 
{Logarithmic bundles and line arrangements, an approach via the standard construction}, \textit{J. London.Math.Soc.}
\textbf{90} (2014), {675--694}.


\bibitem{MuS}  M.\ Musta\c t\u a, H.\ Schenck,
The module of logarithmic p-forms of a locally free arrangement,
\textit{J.\ Algebra} \textbf{241} (2001), 699--719.


\bibitem{OSS} C. Okonek, M. Schneider, H. Spindler: \emph{Vector Bundles on Complex Projective Spaces}. Progress in Math. n. 3, Birkhauser (1980).

\bibitem{OT} P. Orlik, H. Terao :  \emph{Arrangements of
Hyperplanes}. Grundlehren Math. Wiss., \textbf{300},
Springer-Verlag, Berlin  (1992)

\bibitem{S}
K. Saito, 
Theory of logarithmic differential forms and logarithmic vector fields.
\textit{J. Fac. Sci. Univ. Tokyo} \textbf{27} (1980), 265--291.   


\bibitem{Se} E. Sernesi:  The local cohomology of the jacobian ring, \textit{Documenta Mathematica},  \textbf{19} (2014), 541-565. 

\bibitem{T1}
H. Terao, 
Arrangements of hyperplanes and their freeness I, II. 
\textit{J. Fac. Sci. Univ. Tokyo} \textbf{27} (1980), 293--320.   

\bibitem{T2}
H. Terao, 
Generalized exponents of a free arrangement of hyperplanes and
Shephard-Todd-Brieskorn formula. \textit{Invent. Math}. 
\textbf{63}  (1981),
159--179.

\bibitem{W}
A. Wakamiko, 
On the exponents of 2-multiarrangements. 
\textit{Tokyo J. Math}. \textbf{30} (2007), no. 1, 99--116. 

\bibitem{WY}
M. Wakefield and S. Yuzvinsky,
Derivations of an effective divisor on the complex projective line.
\textit{Trans. Amer. Math. Soc}. \textbf{359} (2007), 4389--4403. 

\bibitem{Y1}
M. Yoshinaga,
Characterization of a free arrangement and
conjecture of
Edelman and Reiner. \textit{Invent. Math.} \textbf{157} (2004), no. 2,
449--454.

\bibitem{Y2}
M. Yoshinaga, 
On the freeness of 3-arrangements. 
\textit{Bull. London Math. Soc.} \textbf{37} (2005), no. 1, 126--134. 



\bibitem{Y3}
M. Yoshinaga, 
Freeness of hyperplane arrangements and related topics. 
\textit{Annales de la Facult\`{e} des Sciences de Toulouse } \textbf{23} (2014), no. 2, 483--512. 

\bibitem{Z}
G. M. Ziegler, 
Multiarrangements of hyperplanes and their freeness.  Singularities (Iowa City, IA, 1986),  345--359,
Contemp. Math., {\bf 90}, Amer. Math. Soc., Providence, RI, 1989. 


\end{thebibliography}
\end{document}